  \documentclass{gretsi}
 \usepackage[latin1]{inputenc}   %
 \usepackage[english,french]{babel}   %
  \usepackage{times}			%
 
\usepackage{subfig}  
\usepackage{graphicx}
\usepackage{color}       %

\usepackage{verbatim}
\usepackage{amsmath}
\usepackage{amsfonts}
\usepackage{url}
\def\NLM{\mbox{NL-means}}
\def\SFNLM{\mbox{SFNL-means}}
\def\DCT{\mbox{DCT-means}}

\def\s2{\sigma^2}
\newcommand\R{\mathbb R}

\titre{Un algorithme de débruitage {\em Non-Local means} espace-fréquence}

\auteur{\coord{Simon}{Postec}{},
        \coord{Jacques}{Froment}{},
    \coord{Béatrice}{Vedel}{}}

\adresse{
\affil{}{Université de Bretagne-Sud, UMR CNRS 6205,
Laboratoire de Mathématiques de Bretagne Atlantique, \\ Campus de Tohannic BP 573, F-56017 Vannes}
}

\email{Prenom.Nom@univ-ubs.fr}

\resumefrancais{Le principe de l'algorithme de débruitage {\em Non-Local means} ($\NLM$) \cite{bcm2005, bcm2010} repose sur l'identification dans
l'image bruitée de pixels similaires à partir de la détection de patchs eux-mêmes similaires. Ainsi les détails fins et les motifs 
peu contrastés sont mal restaurés par $\NLM$ mais, comme ils tendent à générer des structures redondantes dans le domaine de Fourier,
un filtrage  $\NLM$ dans ce domaine permet de mieux les débruiter. Une approche mixte espace-fréquence améliore ainsi
les performances globales de $\NLM$, en assurant que toute l'information soit structurée de manière suffisamment redondante, dans le domaine spatial ou dans le domaine fréquentiel. Des approches basées sur l'application de $\NLM$  dans le domaine des ondelettes \cite{sba2006, sn2011} ou sur 
l'ajout d'une contrainte de similarité dans le domaine de Fourier pour le débruitage par variation totale non-locale \cite{hf2012} ont déjà été proposées. Cependant l'article \cite{sba2006} ne fournit aucun résultat quantitatif permettant d'évaluer ses performances et nos résultats se comparent favorablement à ceux présentés dans \cite{sn2011} et \cite{hf2012}. L'intérêt de notre approche réside également dans la simplicité de l'algorithme, qui consiste à appliquer deux fois $\NLM$ (dans le domaine fréquentiel puis dans le domaine spatial). Par rapport à une simple application du filtre $\NLM$ dans le domaine spatial, nous obtenons une meilleure restauration visuelle des textures fines ainsi que des points isolés et le gain global en rapport signal à bruit de crête (PSNR) peut
dépasser 1 dB. Cette approche espace-fréquence procure des résultats intermédiaires entre l'algorithme $\NLM$ original et les méthodes qui font 
l'état de l'art, tout en restant d'une complexité raisonnable et en présentant peu d'artefacts visuels.}

\resumeanglais{{\em A space-frequency Non-Local means image denoising algorithm}. The efficiency of the {\em Non-Local means} ($\NLM$) \cite{bcm2005, bcm2010} image denoising algorithm relies on the identification of similar original pixels from noisy similar patches. Hence fine details and low-contrasted structures are badly recovered after the application of $\NLM$. But as these structures tend to correspond to redundant ones in the Fourier domain, $\NLM$ filtering in this domain allows one to better denoise them. A mixed space-frequency approach improves the denoising performances of $\NLM$ because it ensures that the information is redundant enough, in the spatial domain or in the frequency domain. Approaches based on filtering with $\NLM$ in the wavelet domain \cite{sba2006, sn2011} or on adding a similarity constraint in the Fourier domain for non-local TV denoising \cite{hf2012} have already been proposed. But in \cite{sba2006} there is no experimental result to assess the denoising performances and our results compare favorably to the ones given in \cite{sn2011} and \cite{hf2012}. Our approach is simple : it consists in running two times the $\NLM$ algorithm (firstly in the frequency domain and secondly in the spatial domain). For fine textures and isolated points we get a better visual reconstruction than with the original $\NLM$. In terms of PSNR, the improvement can be over 1 dB. Our approach gives intermediate results between the original $\NLM$ and state-of-the-art methods while at the same time having moderate complexity and leading to few visual artifacts.}

\begin{document}
\maketitle

\section{L'algorithme {\em Non-Local means} en espace et en fréquence}

Soient $v : \Omega \longrightarrow \R$ l'image bruitée et $u$ l'image originale. 
Nous supposons que $v = u + b$, où  $\{b(x)\}_{x \ \in \ \Omega}$ est un bruit blanc gaussien de variance $\s2$.
Le débruitage spatial par l'algorithme {\em Non-Local means} ($\NLM$) \cite{bcm2005} s'écrit classiquement
\begin{equation}
\label{eq:nlm}
\forall x \in \Omega, \ \,\; 
\mbox{NLM}(v)(x) = \displaystyle{\sum_{y \ \in \ \Omega_{x} } \frac{1}{Z(x)} e^{-\frac{\|V(x) - V(y)\|_{2,a}^{2}}{2{h}^2}} v(y)},
\end{equation}
où $V(x)$ est le patch de taille $7 \times 7$ centré en $x$ dans l'image $v$ c'est-à-dire les niveaux de gris $\{v(y) \ | \ y\in\Omega,{\| x - y \|}_\infty \leq 3\}$, $h$ un paramètre de filtrage, 
$Z(x)$ un coefficient de normalisation, ${\| \ . \ \|}_{2,a}$ la norme euclidienne pondérée par une gaussienne d'écart-type 
$a$ qui définit la similarité entre les patchs et $\Omega_x = \{ \ y  \in  \Omega \ | \ {\| x - y \|}_2 \leq d \}$ 
est la zone de recherche des similarités, de rayon $d$.

\begin{figure*}[htb]
\begin{center}
\centerline{\includegraphics[width=2\columnwidth]{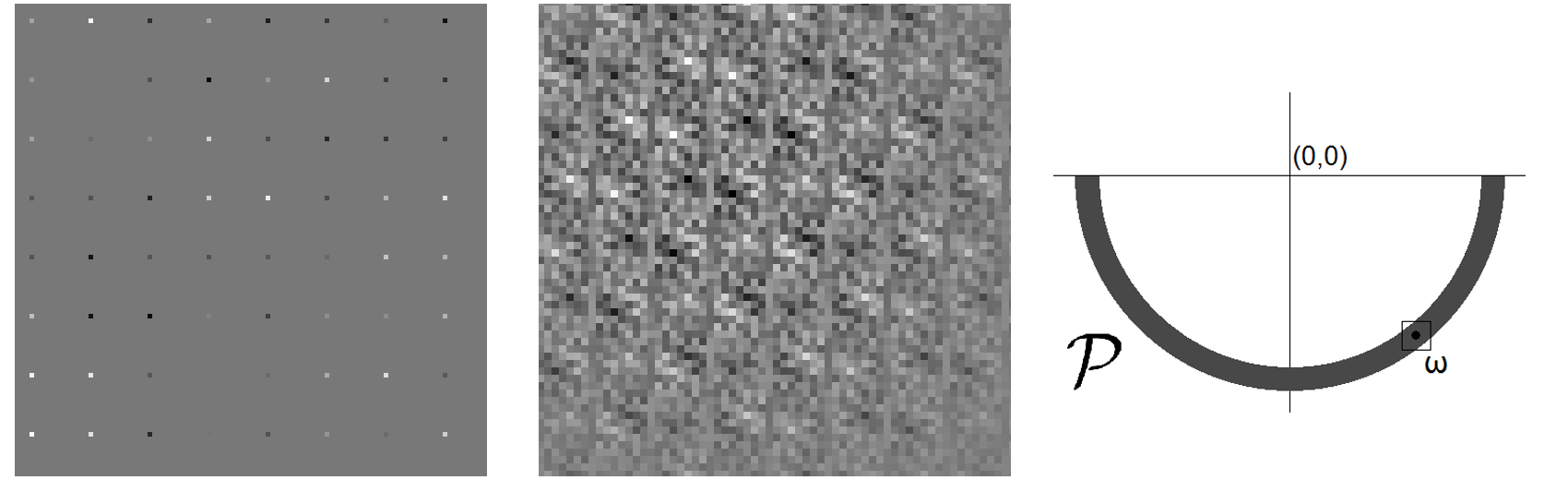}}
\end{center}
\legende{De gauche à droite : image parcimonieuse $p$, $\Re(\hat{p})$ (zoom sur une partie de $\mathcal{P}$), demi-couronne fréquentielle $C_{\omega}$ .}
\label{fig:figure1}
\end{figure*}

Pour être efficace, ce calcul de moyenne pondérée implique de trouver un grand nombre de patchs similaires,
mais certaines structures peuvent appartenir à des patchs qui ne sont pas redondants dans leur voisinage direct
et, comme les meilleures performances de $\NLM$ s'obtiennent en limitant la recherche des patchs similaires à un petit 
voisinage (\cite{p2012}, \cite{pfv2013}), une redondance éventuelle à longue distance n'est en pratique jamais exploitée.
Remarquons qu'aux structures parcimonieuses dans le domaine spatial tendent à correspondre des structures redondantes dans le 
domaine de Fourier, comme cela est illustré à la figure \ref{fig:figure1} : l'image parcimonieuse $p$
n'est pas adaptée à l'algorithme $\NLM$ (bien que présentant de nombreux patchs similaires, l'application de $\NLM$ conduirait 
à un filtrage trop important des points singuliers et donc à la disparition d'une information visuellement significative), 
alors qu'au contraire sa transformée de Fourier $\hat{p}$ donnée à la figure \ref{fig:figure1} contient des structures redondantes bien adaptées à $\NLM$.
Nous proposons donc d'estimer $u$ en introduisant le filtrage $\NLM$ dans le domaine de Fourier suivant le schéma
\begin{equation}
\label{eq:fnlm}
\displaystyle{\mathcal{F}^{-1}(\mbox{FNLM}(\mathcal{F}(v)))}
\end{equation}
où $\mathcal{F}$ désigne la transformée de Fourier 2D discrète, $\mathcal{F}^{-1}$ la transformée inverse et où
le filtrage de l'image $\hat{v}=\mathcal{F}(v)$ à la fréquence $\omega$ s'écrit :
\begin{equation}
\mbox{FNLM}(\hat{v})(\omega) = \displaystyle{\sum_{\xi \ \in C_{\omega} } \frac{1}{Z(\omega)} e^{-\frac{\|\hat{V}(\omega) - \hat{V}(\xi) 
\|_{2,a}^{2}}{2 l^{2}}} \hat{v}(\xi)}.
\label{eq:fnlm_filter}
\end{equation}
Dans cette expression, le filtrage n'est appliqué que dans le demi-plan inférieur du domaine fréquentiel, noté $\mathcal{P}$, 
la reconstruction de l'estimation de $\hat{u}=\mathcal{F}(u)$ s'effectuant grâce aux propriétés de symétrie de la transformée de 
Fourier d'une image réelle (une alternative équivalente serait de remplacer la transformée de Fourier par celle en cosinus).
Pour $\omega \in \mathcal{P}$, $\hat{V}(\omega)$ est le patch de taille $7 \times 7$ centré en la fréquence $\omega$ dans l'image 
$\hat{v}$, $l$ est un paramètre de filtrage, $Z(\omega)$ est un coefficient de normalisation, 
$\|\hat{V}(\omega) - \hat{V}(\xi) \|_{2,a}^{2} = \|\Re(\hat{V}(\omega))- \Re(\hat{V}(\xi))\|_{2,a}^{2} + 
\|\Im(\hat{V}(\omega))- \Im(\hat{V}(\xi))\|_{2,a}^{2}$ où $\Re$ et $\Im$ sont les parties réelles et imaginaires et 
$C_{\omega} = \{ \xi \in \mathcal{P} \ | \ | |\omega| - |\xi| | \leq r \}$ est la demi-couronne de rayon $r$
en la fréquence $\omega$, cf. figure \ref{fig:figure1}.

\section{L'algorithme {\em Space-Frequency\\ Non-Local Means} ($\SFNLM$)}

\begin{figure}[htb]
\begin{center}
\includegraphics[width=85mm]{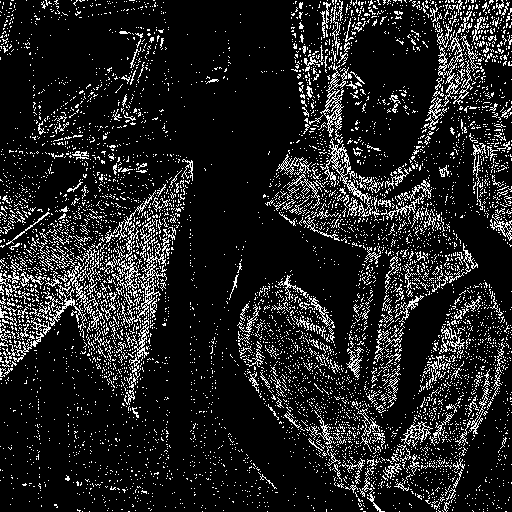}
\end{center}
\legende{Carte, sur l'image {\em Barbara} et en moyenne sur 10 réalisations du bruit, des pixels mieux débruités en Fourier: un pixel $x$ est affiché en blanc si $\sum_{b} |u - \mbox{NLM}(v)|^2 (x) > \sum_{b} |u - \mathcal{F}^{-1}(\mbox{FNLM}(\hat{v}))|^2(x)$.}
\label{fig:figure2}
\end{figure}

Au sens du rapport signal à bruit de crête (PSNR), le filtrage dans le domaine fréquentiel s'avère préférable dans les zones 
texturées et les fins contours alors que le filtrage spatial donne de meilleurs résultats dans les zones homogènes 
(cf. figure \ref{fig:figure2}). Pour autant et comme mentionné précédemment dans le cas de l'image parcimonieuse de la figure 
\ref{fig:figure1}, dans les zones homogènes le filtrage $\NLM$ classique a tendance à effacer des détails ponctuels 
perceptuellement significatifs. Pour ces raisons, nous proposons d'effectuer l'essentiel du filtrage dans le domaine de Fourier
et, afin de régulariser les zones homogènes, de faire suivre ce premier traitement d'un filtrage plus modéré dans le domaine spatial.

\begin{figure}[htb]
\begin{center}
\includegraphics[width=85mm]{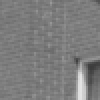}
\end{center}
\legende{Extrait de l'image {\em House}, dont le mur présente une texture difficile à débruiter.\\ \\ \\}
\label{fig:figure3}
\end{figure}

Nous appelons l'algorithme correspondant \textit{Space-Frequency Non-Local Means} ($\mbox{\SFNLM}$) et il s'écrit formellement
\begin{equation}
\label{eq:sfnlm}
\mbox{SFNLM} (v) = \mbox{NLM}(\displaystyle{\mathcal{F}^{-1}(\mbox{FNLM}(\mathcal{F}(v)))}).
\end{equation}
Après quelques essais pour obtenir des résultats homogènes sur différentes images, les paramètres ont été fixés aux valeurs suivantes: 
$l=0.8\sigma$, $r=2$ pour le moyennage fréquentiel; $h=0.6\sigma$ et $d=4$ pour le moyennage spatial, 
\`A condition qu'il soit modéré, le filtrage dans le domaine spatial permet de débruiter les zones homogènes de l'image qui sont mal 
restaurées par le filtre fréquentiel, tout en préservant les points isolés et les textures fines, que le filtre fréquentiel a permis 
de débruiter (cf. figures \ref{fig:figure3} à \ref{fig:figure6} où nous comparons, sur un extrait de l'image {\em House} et pour  $\sigma=10$, 
les résultats donnés par les algorithmes $\NLM$ et $\SFNLM$). 
Les résultats présentés dans la table \ref{tab:arraypsnr20} confirment ces bonnes performances au sens du critère objectif du PSNR. 
Précisons que dans les colonnes \cite{hf2012} et \cite{sn2011} les résultats que nous affichons correspondent à ceux reportés dans ces articles.

\section{Conclusion et perspectives}

La simplicité du schéma $\NLM$ explique son succès face à d'autres algorithmes de débruitage, aux performances
supé\-rieures mais dont la complexité à la fois conceptuelle et algorithmique peut être pénalisante \cite{dfke2007, ea2006}. 
Une meilleure compré\-hension des faiblesses de $\NLM$ permet de proposer des variantes plus performantes au prix d'un accroissement
modéré de la complexité. Ainsi, le manque d'adaptation de $\NLM$ aux structures de l'image peut être corrigé en 
optimisant localement les paramètres de filtrage (voir par exemple \cite{kb2006}). Nous proposons une alternative,
basée sur le modèle statistique sous-jacent dans $\NLM$ de redondance des structures. Le principe est de débruiter
les coefficients obtenus par l'application de transformations inversibles complémentaires au sens des structures redondantes qu'elles 
révèlent. Comme $\NLM$ ne peut exploiter cette redondance que dans un petit voisinage (pour une question d'efficacité algorithmique mais aussi parce
que cet algorithme est optimal lorsqu'il est local \cite{pfv2013}), ces transformations devraient être choisies pour leur capacité à générer des redondances locales.
Dans cette communication, nous illustrons ce principe avec deux transformations élémentaires correspondant à une décomposition dans la base 
canonique et dans la base des exponentielles complexes, mais des décomposi\-tions plus sophistiquées pourraient être appliquées.

Notons que ces deux bases (canonique et exponentielles complexes) forment ce que l'on appelle des dictionnaires incohérents, très étudiés récemment, notamment dans le cadre du \textit{compressed sensing} (\cite{crt2006}) et de la théorie de l'incertitude (\cite{rt2012}).
Le fait que la base canonique et la base des exponentielles complexes forment des dictionnaires à l'incohérence maximale explique certainement les bonnes performances de l'algorithme $\SFNLM$. 

Incidemment, remarquons qu'un filtrage local des coefficients de Fourier entraînant une modification globale des coefficients dans la base canonique,
la méthode de débruitage $\SFNLM$ s'avère être un algorithme non local.

\begin{table}
\legende{ \label{tab:arraypsnr20}Comparaison, en termes de PSNR, d'algorithmes basés sur une régularisation dans le domaine spatial et/ou fréquentiel, pour $\sigma=20$. 
Pour l'algorithme $\NLM$, les valeurs des paramètres sont fixées suivant \cite{pfv2013} à $h = \sigma$ et $d=4$.}
\begin{center}
\small\addtolength{\tabcolsep}{-5pt}
\begin{tabular}{|c||c|c|c|c|c|}
\hline
& $\SFNLM$ & $\NLM$ \cite{bcm2005} & \cite{hf2012} & $\DCT$ \cite{ys2011} & \cite{sn2011} \\
\hline
$Lena$ & \textbf{32.2} & 31.6  & 31.8 & 32.1 & 31.6 \\
\hline
$Barbara$ & 30.0 & 29.2  & 29.2 & 29.9 & \textbf{30.3}\\
\hline
$House$ & \textbf{32.7} & 32.1  & 32.1 & 32.3 & \\
\hline
$Mandrill$ & \textbf{25.9} & 25.8  &  & 25.8 &\\
\hline
$Peppers$ & \textbf{30.6} & 30.4  & 30.3 & 30.0 & 29.2\\
\hline
$Cameraman$ & \textbf{29.6} & 29.4  & 29.6 & 29.1 &\\
\hline
\end{tabular}
\end{center}
\end{table}

\begin{figure}[htb]
\begin{center}
\includegraphics[width=85mm]{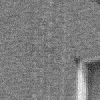}
\end{center}
\legende{Extrait de l'image bruitée, $\sigma=10$, PSNR = 28.14 (image entière).\\ \\ \\}
\label{fig:figure4}
\end{figure}

\begin{figure}[htb]
\begin{center}
\includegraphics[width=85mm]{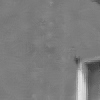}
\end{center}
\legende{Extrait de la restauration par $\NLM$, PSNR = 36.16 (image entière). 
Noter la disparition presque complète de la texture du mur, ce qui donnerait à penser que le paramètre
de lissage $h$ est trop élevé, tandis qu'au contraire le contour de la fenêtre reste bruité.}
\label{fig:figure5}
\end{figure}

\begin{figure}[htb]
\begin{center}
\includegraphics[width=85mm]{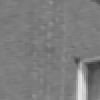}
\end{center}
\legende{Extrait de la restauration par $\SFNLM$, PSNR = 37.19 (image entière).
La texture du mur et le contour de la fenêtre sont mieux préservés. 
Le critère objectif du PSNR prend bien en compte l'amélioration perceptive avec un accroissement de 
plus de 1 dB (pour l'ensemble de l'image).
}
\label{fig:figure6}
\end{figure}

\end{document}